\documentclass[12pt,leqno]{article}
\usepackage{amssymb}
\usepackage[mathscr]{eucal}
\usepackage{amsmath,amssymb,latexsym,theorem,bbm}
\usepackage{color,url}

\usepackage{mathtools}
\mathtoolsset{mathic=true}
\usepackage{microtype}

\newcommand{\NN}{\mathbb{N}}

\newcommand{\RR}{\mathbb{R}}

\newcommand{\ZZ}{\mathbb{Z}}

\newcommand{\bA}{{\boldsymbol{A}}}

\newcommand{\cA}{{\mathcal A}}

\newcommand{\cD}{{\mathcal D}}

\newcommand{\dd}{\mathrm{d}}
\newcommand{\ee}{\mathrm{e}}

\newcommand{\PP}{\operatorname{\mathbb{P}}}

\newcommand{\tZ}{\widetilde{Z}}

\newcommand{\vare}{\varepsilon}

\renewcommand{\leq}{\leqslant}
\renewcommand{\geq}{\geqslant}

\newcommand{\distre}{\stackrel{\cD}{=}}

\newcommand{\bbone}{\mathbbm{1}}

\newcommand{\proofend}{\hfill\mbox{$\Box$}}

\numberwithin{equation}{section}

\theoremstyle{change} \theorembodyfont{\em}
\newtheorem{Lem}{Lemma.}[section]
\newtheorem{Thm}[Lem]{Theorem.}

\newtheorem{Cor}[Lem]{Corollary.}

\theorembodyfont{\rm}
\newtheorem{Rem}[Lem]{Remark.}

\begin{document}

\begin{center}
{\bfseries\Large Pick-up sticks problems\\[1mm] for exponential and Pareto distributions}
\\[1mm]
{\bfseries\Large }
\vspace*{3mm}

 {\sc\large
  M\'aty\'as $\text{Barczy}^{*,\diamond}$, Rezs\H o L. $\text{Lovas}^{**}$ }

\end{center}

\vskip0.2cm

\noindent
 * HUN-REN--SZTE Analysis and Applications Research Group,
   Bolyai Institute, University of Szeged,
   Aradi v\'ertan\'uk tere 1, H--6720 Szeged, Hungary.

\noindent
 ** Institute of Mathematics, University of Debrecen,
    Pf.~400, H--4002 Debrecen, Hungary.

\noindent E-mails: barczy@math.u-szeged.hu (M. Barczy),
                   lovas@science.unideb.hu (R. L. Lovas).

\noindent $\diamond$ Corresponding author.

\renewcommand{\thefootnote}{}

\vspace*{-0.5cm}

\begin{abstract}
We deal with the following two pick-up sticks problems:
 given an integer $n$ greater than or equal to $3$, we choose $n$ sticks with independent 
 random lengths from an exponential distribution or Pareto type I distribution,
 and we are interested in the probability that no three of them can be the sides of a triangle.
In case of exponential distribution, this problem has been solved by other authors:
 the probability in question is equal to $n!/\big((F_3-1)(F_4-1)\cdots(F_{n+2}-1)\big)$, 
 where $(F_i)_{i=1}^\infty$ is the sequence of Fibonacci numbers.
We give a new, self-contained and detailed proof for this formula, which uses a different method compared to the earlier ones.
In case of Pareto type I distribution, we derive new formulae for the corresponding probability.
\end{abstract}

\section{Introduction}

Pick-up sticks problems are the natural counterparts of the well-known and well-studied broken stick problems.
According to our knowledge, the pick-up sticks model has been recently introduced by 
 Petersen and Tenner \cite{PetTen}, who originally considered the following problem:
 given an integer $n$ greater than or equal to $3$, we choose $n$ random numbers independently 
 according to uniform distribution on $[0,1]$ (there is no constraint that they should sum to one),
 and we are interested in the probability that $n$ such sticks cannot form an $n$-gon
 (by an $n$-gon, we mean a (convex or concave) polygon with $n$ sides and $n$ vertices
 whose sides do not intersect with each other).
A well-known generalized triangle inequality for polygons says that 
 a collection of positive numbers $\ell_1,\ldots,\ell_n$ can be the side lengths of an $n$-gon if and only if
 $\ell_i < \sum_{j\in\{1,\ldots,n\}\setminus\{i\}} \ell_j$, $i=1,\ldots,n$, or equivalently, 
 $\max_{i\in\{1,\ldots,n\}}\ell_i < \frac{1}{2}(\ell_1+\cdots+\ell_n)$.
 
In this paper, we deal with the following two variants of the pick-up sticks problem recalled above:
 given an integer $n$ greater than or equal to $3$, we choose $n$ sticks with independent 
 random lengths from an exponential distribution or Pareto type I distribution, and we are interested in
 the probability that no three of them can be the lengths of a triangle.
Below, we refer to and summarize only the relevant literature that are connected to our problems.
There is a much more extensive body of literature available on the broken sticks problems.

First, consider the following broken stick problem:
 we break a stick at random $n-1$ points to form $n$ pieces, and we are interested in the probability that no choice
 of $k$ pieces out of $n$ forms a $k$-gon, where $n$ is an integer greater than or equal to 3, and $k\in\{3,\ldots,n\}$.
This problem was solved by Verreault \cite[Theorem 2]{Ver2} using combinatorial techniques
 based on MacMahon partition analysis, 
 and, later on, by Mukerjee \cite[Theorem 1]{Muk},
 through a direct statistical approach.
For historical fidelity, we mention that a formula for the probability in question was stated earlier \
 without any proof in the research project report by Kong et al.\ \cite[Theorem 1]{KonLkhTurUppHil} of the Illinois Geometry Lab.
The proof of Mukerjee \cite[Theorem 1]{Muk} is based on the following fact: 
 the order statistics corresponding to the $n$ lengths of the $n$ pieces when a stick having length
 1 is broken at random at $n-1$ points (chosen according to uniform distribution on $[0,1]$)
 has the same joint distribution as 
  \begin{align}\label{help_exp_4}
    \left(\frac{Z_{(1)}}{\sum_{j=1}^n Z_j},\ldots,\frac{Z_{(n)}}{\sum_{j=1}^n Z_j}\right),
  \end{align}
 where $Z_1,\ldots,Z_n$ are independent and exponentially distributed random variables with parameter $1$,
 and $Z_{(1)},\ldots,Z_{(n)}$ satisfying $Z_{(1)}\leq Z_{(2)}\leq \cdots\leq Z_{(n)}$ denote their order statistics, see, e.g., 
 Verreault \cite[Lemma 1]{Ver1} or Mukerjee \cite[Proposition 1]{Muk}.
Then, using the well-known formula for the joint density function of the order statistics for exponentially distributed random 
 variables, Mukerjee \cite{Muk} expressed the probability in question as an integral, where a linear variable transformation
 was carried out to obtain a closed form expression.
In the Appendix of Mukerjee \cite{Muk}, he also checked that his formula in Theorem 1 indeed coincides with the formula derived
 earlier by Verreault \cite[Theorem 2]{Ver2}, where generalized Fibonacci numbers came into play.

The reason for presenting the previous broken stick problem and the methods for its solution in a detailed way is that 
 it is in fact equivalent to the pick-up sticks problem due to Sudbury et al.\ \cite{SudSunTreWan} (presented below)
 when the uniform distribution is replaced by exponential distribution. 
This follows from the fact that the order statistics corresponding to the $n$ lengths of the broken 
 stick problem in question can be represented as in \eqref{help_exp_4}, and we mention that
 the equivalence of the two problems in question was also noted in Brennan et al.\ \cite[page 19]{BreCalLin}. 

Before recalling the pick-up sticks problem considered by Sudburry et al.\ \cite{SudSunTreWan}, we mention that
Liu and Wang \cite[Theorem 2.4]{LiuWan} investigated and solved the following pick-up sticks problem:
 given an integer $n$ greater than or equal to $3$, if we choose independently $n$ sticks with lengths from 
 uniform distribution on $[0,1]$, then the probability that any three sticks among them can form a triangle
 is $1/2^{n-2}$.

Sudbury et al.\ \cite{SudSunTreWan} considered the following variant of the pick-up sticks model:
 given an integer $n$ greater than or equal to $3$, if we choose $n$ sticks with independent random lengths
 selected from a uniform distribution on $[0,1]$, what is the probability that no three of them can be the sides of a triangle?
It was calculated that this probability is $1/(F_1\cdots F_n)$, where $(F_i)_{i=1}^\infty$ is the sequence of Fibonacci numbers
 defined by $F_1:=1$, $F_2:=1$, and $F_i := F_{i-1}+ F_{i-2}$, $i\geq 3$. 
Concerning the method of Sudbury et al.\ \cite{SudSunTreWan}, they used that
 the event that no three of the selected sticks can be the sides of a triangle is equivalent that all the consecutive triplets
 $(\ell_{(i-2)},\ell_{(i-1)},\ell_{(i)})$ satisfy the inequality $\ell_{(i-2)}+\ell_{(i-1)}\leq \ell_{(i)}$, $i=3,\ldots,n$,
 where $0\leq \ell_{(1)}\leq \ell_{(2)}\leq \cdots\leq \ell_{(n)}$ denotes the increasing order of lengths.
Then they used the well-known fact that the joint distribution of 
 the order statistics of an $n$-sample from a uniform distribution on $[0,1]$ 
 coincides with that of  the normalized cumulative sums 
 of $n+1$ independent exponentially distributed random variables with parameter 1
 (see, e.g., Shorack and Wellner \cite[Chapter 8, Proposition 1]{ShoWel}),
 i.e., if $n$ is a positive integer, $U_1,\ldots,U_n$ are independent and identically distributed random variables 
 such that their common distribution is uniform on $[0,1]$, and $Z_1,\ldots,Z_{n+1}$ are independent and identically distributed random
 variables such that their common distribution is exponential with parameter $1$, then
 $(U_{(1)},\ldots, U_{(n)})$ has the same distribution as
 \[
   \left( \frac{Z_1}{\sum_{i=1}^{n+1}Z_i}, \, \frac{Z_1 + Z_2}{\sum_{i=1}^{n+1} Z_i},\, \ldots, \, \frac{Z_1 + \dots + Z_n}{\sum_{i=1}^{n+1} Z_i} \right),
 \]
 where $U_{(1)},\ldots,U_{(n)}$ satisfying $U_{(1)}\leq U_{(2)}\leq \cdots\leq U_{(n)}$ denote the order statistics of $U_1,\ldots,U_n$.
Using these observations, Sudbury et al.\ \cite[Theorem 1]{SudSunTreWan}  expressed the probability in question as a multiple integral, 
 and, by successively evaluating it, they arrived at the 
 reciprocal of the product of the first $n$ Fibonacci numbers.
Sudbury et al.\ \cite[Theorem 2]{SudSunTreWan} also considered the case of quadrilaterals instead of triangles, 
 and derived a formula for the corresponding probability in terms of 
 Tribonacci (3-step Fibonacci) numbers.

Brennan et al.\ \cite[Theorem 1.1]{BreCalLin} have recently generalized Theorems 1 and 2 in  Sudbury et al.\ \cite{SudSunTreWan}.
Namely, given $n$ sticks, which have independent random lengths chosen from a uniform distribution on $[0,1]$,
 they calculated the probability that no $k$ of them can form a $k$-gon, where $n$ is a positive integer greater than or equal to $3$ and $k\in\{3,\ldots,n\}$.
The probability in question depends on the sequence $(F^{(k-1)}_i)_{i=-k+3}^\infty$ of $(k-1)$-step Fibonacci numbers defined by the recursion 
 $F^{(k-1)}_0:=F^{(k-1)}_{-1}:=\cdots = F^{(k-1)}_{-k+3}:=0$, $F^{(k-1)}_1:=1$, and 
 $F^{(k-1)}_i:=\sum_{j=1}^{k-1} F^{(k-1)}_{i-j}$, $i\geq 2$. 
Note that, in case of $k=3$, $(F^{(2)}_i)_{i=1}^\infty$ is nothing else but the sequence $(F_i)_{i=1}^\infty$ of (usual) Fibonacci numbers.
Brennan et al.\ \cite[Theorem 1.1]{BreCalLin} not only generalized the results in Sudbury et al.\ \cite{SudSunTreWan}, but they also developed an alternative
 geometric proof.
In case of $k=3$ (triangles), the essence of their method can be summarized as follows.
Let us sort the $n$ sticks in order of increasing lengths $\ell_{(1)}\leq \ell_{(2)} \leq \cdots \leq \ell_{(n)}$.
Then, as it was explained earlier, the event that no three of the selected sticks can be the sides of a triangle is equivalent that
 the inequalities $\ell_{(i)} \geq \ell_{(i-2)} + \ell_{(i-1)}$, $i=3,\ldots,n$, hold,
 which means that there is a lower bound for $\ell_{(i)}$, namely, $\ell_{(i-2)} + \ell_{(i-1)}$.
The key observation in the proof of Theorem 1.1 in  Brennan et al.\ \cite{BreCalLin} (in case of triangles) 
 is that the above lower bounds for $\ell_{(i)}$, $i=3,\ldots,n$, also induce constraints on the maximum lengths of the sticks,
 namely, $\ell_{(n)}\leq 1$ and 
 \[
  \ell_{(i)}\leq \frac{1-F_{n-i}\ell_{(i-1)}}{F_{n+1-i}},\qquad i=1,\ldots,n-1,
 \]
 where $\ell_{(0)}:=0$, see Brennan et al.\ \cite[formula (3.3)]{BreCalLin}. 
They also pointed out the fact that the inequality above explains that Fibonacci numbers appear naturally in Theorem 1 in  
 Sudbury et al.\ \cite{SudSunTreWan}.
The same method was applied by Brennan et al.\ \cite[formula (3.9) and the proof of Lemma 3.6]{BreCalLin} 
 in case of $k$-gons as well (where $k\in\{3,\ldots,n\}$).

Brennan et al.\ \cite[formula (3.24)]{BreCalLin} also derived a formula for the probability 
 that no $k$ of $n$ sticks, which have independent random lengths chosen from an unbounded nonnegative absolutely continuous 
 distribution, can form a $k$-gon, where $n$ is an integer greater than or equal to $3$ and $k\in\{3,\ldots,n\}$.
This formula is in fact an integral, which depends on the density function of the absolutely continuous distribution
 in question.
As particular examples, they considered the cases of exponential distribution and shifted exponential distribution
 with support $[b,\infty)$, where $b>0$, and in these cases explicit formulae were possible to obtain,
 see Brennan et al.\ \cite[formula (3.26)]{BreCalLin}.
In case of exponential distribution, the formula depends on the sequence of $(k-1)$-step Fibonacci numbers (recalled earlier), 
 which, in case  of triangles ($k=3$), reduces to $n!/\big((F_3-1)(F_4-1)\cdots(F_{n+2}-1)\big)$, where $(F_i)_{i=1}^\infty$ is the sequence of
 (usual) Fibonacci numbers, see page 20 in Brennan et al.\ \cite{BreCalLin}.
In case of exponential distribution and $k=3$, the method of evaluating the integral for the corresponding probability is the one 
 with minor modifications that we previously recalled for the pick-up sticks problem for uniform distribution (lower and upper bounds for the lengths of
 the sticks that are compatible with the requirement that no three of them can be the side lengths of a triangle),
 see the proof of Theorem 4.4 in Brennan et al.\ \cite{BreCalLin}.
For completeness, we mention that Theorem 4.4 in Brennan et al.\ \cite{BreCalLin} has been formulated for the following broken sticks problem:
 a stick of unit length is broken at $(n-1)$ random points (according to uniform distribution on $[0,1]$) to form $n$ random length pieces, and we are interested in
 the probability of not being able to form a $k$-gon from any possible combination of $k$ of these pieces, where $k\in\{3,\ldots,n\}$.
As we already pointed out, as a consequence of \eqref{help_exp_4}, in case of $k=3$ this broken sticks problem is equivalent 
 to the pick-up sticks problem in question (choosing $n$ sticks with independent random lengths selected from exponential distribution,
 what is the probability that no three of them can be the side lengths of a triangle?).

The paper is structured as follows.
In Section \ref{Sec_exp}, we give a new, self-contained and detailed proof of the following pick-up sticks problem
 considered earlier by other authors (detailed in the previous paragraphs):
 given $n$ sticks, which have independent random lengths chosen from an exponential distribution,
 we calculate the probability that no three of them can form a triangle, see Theorem \ref{thm2}.
According to our knowledge, our proof is different from all the existing ones in the literature, namely, 
 from the ones in Sudbury et al.\ \cite{SudSunTreWan} (who considered the pick-up sticks problem in question for uniform distribution), 
 in Brennan et al.\ \cite{BreCalLin} (who considered the pick-up sticks problem in question for exponential distribution),
 in  Mukerjee \cite{Muk} and in Verreault \cite{Ver2} (who considered the equivalent broken sticks problem). 
Our proof of Theorem \ref{thm2} is based on R\'enyi's representation of order statistics for exponential
 distribution (and not for uniform distribution) in terms of some independent exponentially distributed random variables,
 see \eqref{Renyi_rep}.
Using R\'enyi's representation and the joint density function of the independent exponentially distributed random variables in question,
 we write the probability in question as an integral, see \eqref{help_exp_3}.
This integral is evaluated successively, but without performing any change of variables.
In the course of successive integration, a matrix $A_n$ (defined in Step 9 in the proof of Theorem \ref{thm2})
 plays a crucial role, and its form sheds some light on the fact that the sums $F_1+\cdots+F_i = F_{i+2}-1$, $i=1,\ldots,n$,
 of Fibonacci numbers appear in the explicit expression of the corresponding probability. 
As we mentioned earlier,  Brennan et al.\ \cite[formula (3.3)]{BreCalLin} also gave another explanation for the natural 
 appearance of Fibonacci numbers.
 
In Section \ref{Sec_Par}, we consider a variant of the pick-up sticks problem investigated in Section \ref{Sec_exp} by replacing the exponential distribution 
 with a Pareto type I distribution. 
Theorem \ref{thm4} contains two formulae (see \eqref{Pareto_formula_3}) for the probability in question in terms of integrals, and our proof is based on
 a representation of the ratios of successive order statistics for Pareto type I distribution with scale parameter 1 and shape parameter $\alpha>0$ (see part (iii) of Lemma \ref{Lem_Pareto}). 
In Corollary \ref{cor}, we specialize Theorem \ref{thm4} to the case $n=3$, and, in addition, when $\alpha$ is a positive integer,
  by evaluating the integrals mentioned before, we derive an explicit formula for the probability in question.
Remark \ref{Rem1} is devoted to present some other possible approaches for determining the probability in Corollary \ref{cor}.

Finally, we introduce the notations that will be used throughout the paper.
Let $\ZZ_+$, $\NN$, $\RR$, $\RR_+$ and $\RR_{++}$ denote the
 set of non-negative integers, positive integers, real numbers, non-negative real numbers and positive real numbers, respectively.
All the random variables will be defined on a common probability space $(\Omega,\cA,\PP)$.
Equality in distribution of random variables is denoted by $\distre$.
Given an $n\in\NN$, and independent and identically distributed random variables $X_1,\ldots, X_n$,
 their order statistics is denoted by $X_{(1)},\ldots,X_{(n)}$ satisfying $X_{(1)}\leq X_{(2)}\leq \cdots \leq X_{(n)}$.

\section{Pick-up sticks for exponential distribution}\label{Sec_exp}

In this section, we reconsider a pick-up sticks problem for exponential distribution studied by other authors (for details, see the Introduction),
 and for the corresponding probability, we give a new, self-contained and detailed proof different from the existing ones.

\begin{Thm}\label{thm2}
Let $n\in\NN$ be an integer greater than or equal to $3$.
If $n$ numbers are chosen independently according to an exponential distribution with
 parameter $\lambda>0$ from the interval $(0,\infty)$, 
 then the probability that no three of the numbers can be the side lengths of a triangle is
 \[
    \prod_{i=1}^{n} \frac{i}{F_{i+2}-1}  = \frac{n!}{(F_{3}-1)(F_{4}-1) \cdots (F_{n+2}-1)},
 \] 
 where $F_n$, $n\in\NN$, are the Fibonacci numbers defined by $F_1 := 1$, $F_2 := 1$, and $F_{n} := F_{n-1} + F_{n-2}$, $n > 2$, $n\in\NN$.
\end{Thm}

\noindent{\bf Proof.}
We divide the proof into several steps.

{\sl Step 1.} 
Let $X_1,\ldots,X_n$ be independent, exponentially distributed random variables with 
 parameter $\lambda>0$.
Recall that $X_{(1)}, X_{(2)}, \ldots, X_{(n)}$ satisfying $X_{(1)} \leq X_{(2)} \leq \ldots \leq X_{(n)}$ 
 denote the order statistics of $X_1,\ldots,X_n$.
Since the exponential distribution with parameter $\lambda$ is absolutely continuous, the probability that 
 $X_{(i)} = X_{(j)}$ for some $i<j$, $i,j\in\{1,\ldots,n\}$, is zero, and therefore, 
 without loss of generality, one can assume that 
 \[
    X_{(1)} < X_{(2)} < \ldots < X_{(n)}.
 \] 
Recall also that if $\lambda>0$ and $Z$ is an exponentially distributed random variable with parameter $1$, then 
 $\frac{1}{\lambda} Z$ is an exponentially distributed random variable with parameter $\lambda$. 
  
{\sl Step 2.}   
We check that the requirement that no three of $X_1, X_2, \ldots, X_n$ 
 can be the side lengths of a triangle is equivalent to
  \begin{equation}\label{eq0_5_exp}
    X_{(a)} + X_{(b)} \leq X_{(c)}
  \end{equation}
 for all indices $a,b,c\in\{1,\ldots,n\}$ satisfying $a < b < c$. 
 
First, suppose that the requirement in question holds, i.e., 
 one cannot choose three of $X_1,\ldots,X_n$ such that the sum 
 of any two of them is strictly greater than the third one.
If $a,b,c\in\{1,\ldots,n\}$ are such that $a < b < c$, then $X_{(a)} < X_{(b)} < X_{(c)}$, yielding that 
 $X_{(a)} + X_{(c)} \geq X_{(c)} > X_{(b)}$ and $X_{(b)} + X_{(c)} \geq X_{(b)} > X_{(a)}$,
 and hence, using also the assumption, the inequality $X_{(a)} + X_{(b)} > X_{(c)}$ cannot hold,
 equivalently \eqref{eq0_5_exp} holds, as desired.
 
Conversely, suppose that \eqref{eq0_5_exp} holds for all indices $a,b,c\in\{1,\ldots,n\}$ with $a < b < c$. 
Then it is not true that the sum of any two of $X_{(a)}, X_{(b)}, X_{(c)}$ is strictly greater than the third one.
Consequently, one cannot choose three of $X_1,\ldots,X_n$ such that the sum of any two of them is strictly greater than the third one, as desired.
 
{\sl Step 3.}  
We check that the inequality \eqref{eq0_5_exp} holds for all indices $a,b,c\in\{1,\ldots,n\}$ satisfying $a < b < c$
 if and only if
 \begin{equation}\label{eq1_exp}
    X_{(i)} + X_{(i+1)} \leq X_{(i+2)}
 \end{equation}
 holds for all $i\in\{1,\ldots, n-2\}$.

If \eqref{eq0_5_exp} holds for all indices $a,b,c\in\{1,\ldots,n\}$ satisfying $a < b < c$, then 
 it readily holds for $a=i$, $b=i+1$ and $c=i+2$, where $i\in\{1,\ldots, n-2\}$, yielding \eqref{eq1_exp}.

Suppose now that \eqref{eq1_exp} holds for all $i\in\{1,\ldots, n-2\}$, and let $a,b,c\in\{1,\ldots,n\}$ 
 satisfying $a<b<c$.
Then, by \eqref{eq1_exp} with $i=c-2$, we have $X_{(c)}\geq X_{(c-1)}+ X_{(c-2)}$, and, since $c-1\geq b$ and $c-2\geq a$, 
 we have $X_{(c-1)} \geq X_{(b)}$ and $X_{(c-2)} \geq X_{(a)}$.
Therefore $X_{(c)}\geq X_{(a)} + X_{(b)}$, i.e., \eqref{eq0_5_exp} holds.

{\sl Step 4.}
We recall the so-called R\'enyi's representation of order statistics for
 exponential distribution:
   \begin{align}\label{Renyi_rep}
   \big( X_{(1)},\ldots, X_{(n)}\big) 
     \distre \left( \frac{1}{\lambda} \sum_{j=1}^i \frac{Z_j}{n-j+1} \right)_{i=1}^n,
 \end{align}
 where $Z_1,\ldots,Z_n$ are independent, exponentially distributed 
 random variables with parameter $1$, see R\'enyi \cite{Ren} or Nevzorov \cite[Representation 3.4]{Nev}.
 
{\sl Step 5.}
We study the case $n=3$. 
We could have omitted this step, since in Step 7, when we consider the general case $n\geq 3$, $n\in\NN$, we do not rely on this special case.
But we decided to keep it in order to shed light on our method of proof.
Then \eqref{Renyi_rep} takes the form
 \[
  \big( X_{(1)},X_{(2)}, X_{(3)}\big) \distre
   \frac{1}{\lambda}\left(\frac{Z_1}{3}, \frac{Z_1}{3}+\frac{Z_2}{2}, \frac{Z_1}{3}+\frac{Z_2}{2}+Z_3\right).
 \] 
Hence, by Steps 1--3, the probability in question is
 \begin{align}\label{help_exp_5}
  \begin{split}
  \PP(X_{(1)}+X_{(2)}\leq X_{(3)})
   &=\PP\left(\frac{1}{\lambda}\left(\frac{Z_1}{3}+\frac{Z_1}{3}+\frac{Z_2}{2}\right)
       \leq \frac{1}{\lambda}\left(\frac{Z_1}{3}+\frac{Z_2}{2}+Z_3\right)\right)\\
   & = \PP\left(\frac{Z_1}{3}\leq Z_3\right)
     = \int_{\{ (z_1,z_3)\in\RR_+^2 \,:\, z_1/3\leq z_3 \} }
      f_{Z_1,Z_3}(z_1,z_3)\,\dd z_1\dd z_3\\
    & = \int_0^\infty \left( \int_{z_1/3}^\infty \ee^{-(z_1+z_3)}\,\dd z_3\right)\dd z_1 \\
   & = \int_0^\infty \ee^{-z_1} \ee^{-\frac{z_1}{3}} \,\dd z_1 
     = \int_0^\infty \ee^{-\frac{4}{3}z_1} \,\dd z_1 
     = \frac{3}{4},
   \end{split}  
 \end{align}
 where $f_{Z_1,Z_3}$ denotes the joint density function of $Z_1$ and $Z_3$.
Hence $\PP(X_{(1)}+X_{(2)}\leq X_{(3)})$ coincides with 
 \[
    \frac{3!}{(F_3-1)(F_4-1)(F_5-1)}
      = \frac{6}{(2-1)(3-1)(5-1)}
      = \frac{3}{4},
 \]  
 yielding that the assertion holds for $n=3$.
By \eqref{help_exp_5}, we can also see that $\PP(X_{(1)}+X_{(2)}\leq X_{(3)}) = \PP\left(Z_1\leq 3Z_3\right)$,
 which shows that for the calculation of the probability $\PP(X_{(1)}+X_{(2)}\leq X_{(3)})$, interestingly, we only need two 
 (and not three) independent, exponentially distributed random variables.

{\sl Step 6.}
Let $G_n:=F_1+\cdots+F_n$, $n\in\NN$.
Then 
 \begin{align}\label{G_n_F_n}
    G_n=F_{n+2}-1, \qquad n\in\NN.
 \end{align}
Indeed, in case of $n=1$, we have that $G_1=F_1=1$ and $F_{3}-1=2-1=1$,
 yielding that $G_1=F_{3}-1$.
Suppose that \eqref{G_n_F_n} holds for $1,\ldots,n$.
Then 
 \[
 G_{n+1} = F_1+\cdots+F_{n+1} = G_n + F_{n+1}
         = F_{n+2}-1 + F_{n+1}
         = F_{n+3}-1,
 \]
 as desired.

{\sl Step 7.} We study the general case $n\geq 3$, $n\in\NN$.
Using Steps 1--4, we need to calculate the probability
 \begin{align*}
 P_n&:=\PP\Big(X_{(i)} + X_{(i+1)} \leq X_{(i+2)} \;\; \text{for each $i\in\{1,\ldots,n-2\}$}\Big)\\
   &= \PP\left( (\tZ_1 + \dots + \tZ_i) +  (\tZ_1 + \dots + \tZ_{i+1})
                     \leq \tZ_1 + \dots + \tZ_{i+2} \;\; \text{for each $i\in\{1,\ldots,n-2\}$}
                     \right)\\
   &= \PP\Big( \tZ_1 + \dots + \tZ_{i} \leq \tZ_{i+2} \;\; \text{for each $i\in\{1,\ldots,n-2\}$} \Big),             
 \end{align*}
 where $\tZ_i:=\frac{Z_i}{n-i+1}$, $i=1,\ldots,n$, are independent random variables such that $\tZ_i$ has an exponential distribution with parameter $n-i+1$.
For better readability, we spell out that the inequalities $\tZ_1 + \dots + \tZ_{i} \leq \tZ_{i+2}$, $i\in\{1,\ldots,n-2\}$,
 mean that 
 \begin{align*}
    \tZ_1 \leq \tZ_3, \qquad \tZ_1 + \tZ_2 \leq \tZ_4, \qquad \ldots \qquad  \tZ_1 + \dots + \tZ_{n-2} \leq \tZ_n.
  \end{align*}

Denoting the joint density function of $(\tZ_1,\ldots,\tZ_n)$ by $f_{\tZ_1,\ldots,\tZ_n}$, and using its explicit form
 \[
   f_{\tZ_1,\ldots,\tZ_n}(x_1,\ldots,x_n) = f_{\tZ_1}(x_1)\cdots f_{\tZ_n}(x_n)
                                      = n!\ee^{-(nx_1 + (n-1)x_2+ \cdots+ 2x_{n-1} + x_n)},\qquad x_1,\ldots,x_n\in\RR_+,
 \] 
 we get that 
 \begin{align}\label{help_exp_3}  
    P_n = n!\int_{\big\{ (x_1,\ldots,x_n)\in\RR_+^n \, : \, x_1+\cdots+x_i\leq x_{i+2}, \, i=1,\ldots,n-2 \big\} } 
               \ee^{-\big(nx_1 + (n-1)x_2+ \cdots+ 2x_{n-1} + x_n\big)} \,\dd x_n\ldots\dd x_1.
 \end{align}

{\sl Step 8.}
For each $n\geq 3$, $n\in\NN$, $i\in\{1,\ldots,n-2\}$, and $x_1,x_2,\ldots,x_n\in\RR_+$, let 
 \begin{align}\label{help_exp_1}
  I_0(x_1,\ldots,x_n):= \ee^{-\big(nx_1 + (n-1)x_2+ \cdots+ 2x_{n-1} + x_n\big)},
 \end{align}
 and 
\[
   I_i(x_1,\ldots,x_{n-i})
        :=\int_{x_1 + \dots + x_{n-i-1}}^{\infty} \dotsi \int_{x_1 + \dots + x_{n-2}}^{\infty} 
             \ee^{-\big( nx_1 + (n-1)x_2+ \cdots+ 2x_{n-1} + x_n \big)} \, \dd x_n \ldots \dd x_{n-i+1}.
 \]
Since $\int_a^\infty \ee^{-bx}\,\dd x =\frac{1}{b} \ee^{-ba}$ for all $a>0$ and $b>0$, we have that,
 for each $n\geq 3$, $n\in\NN$ and $i\in\{0,1,\ldots,n-2\}$,
 there exist $\beta_i, \alpha_{i,1}, \alpha_{i,2},\ldots, \alpha_{i,n-i}\in\RR_{++}$ such that 
 \begin{align}\label{help_exp_2}
  I_i(x_1,\ldots,x_{n-i}) = \frac{1}{\beta_i} \ee^{-(\alpha_{i,1} x_1 + \cdots + \alpha_{i,n-i} x_{n-i})},
       \qquad x_1,\ldots,x_{n-i}\in\RR_+.
 \end{align}
Formally, let $\alpha_{i,j}:=0$ for $j\in\{ n-i+1,\ldots,n\}$.
In case of $i=0$, by \eqref{help_exp_1}, we have 
 \[
  \beta_0=1,\qquad \alpha_{0,1}=n,\qquad  \alpha_{0,2}=n-1,\qquad \cdots \qquad \alpha_{0,n}=1.
 \]
Hence, the first row of the matrix $A_n:=(\alpha_{i,j})_{i=0,\ldots,n-2}^{j=1,\ldots,n}\in \RR_+^{(n-1)\times n}$
 takes the form $(n, n-1, \ldots, 2,1)$.
Further, for each $i\in\{0,\ldots, n-3\}$, we get 
 \begin{align*}
   I_{i+1}(x_1,\ldots, x_{n-i-1})
      &= \int_{x_1+\cdots+x_{n-i-2}}^\infty I_i(x_1,\ldots,x_{n-i})\,\dd x_{n-i}\\
      &= \frac{1}{\beta_i} \int_{x_1+\cdots+x_{n-i-2}}^\infty
         \ee^{-(\alpha_{i,1} x_1 + \cdots + \alpha_{i,n-i} x_{n-i})}\,\dd x_{n-i}\\
      &= \frac{1}{\beta_i}    
          \ee^{-(\alpha_{i,1} x_1 + \cdots + \alpha_{i,n-i-1} x_{n-i-1})}
          \frac{1}{\alpha_{i,n-i}}
          \Big[ - \ee^{-\alpha_{i,n-i} x_{n-i}} \Big]_{x_{n-i} = x_1+\cdots + x_{n-i-2}}^{x_{n-i}=\infty}\\
      &= \frac{1}{\beta_i\alpha_{i,n-i}}    
          \ee^{-(\alpha_{i,1} x_1 + \cdots + \alpha_{i,n-i-1} x_{n-i-1})}  
          \ee^{-\alpha_{i,n-i}(x_1+\cdots + x_{n-i-2})} \\
      & = \frac{1}{\beta_i\alpha_{i,n-i}}        
          \ee^{-\big( (\alpha_{i,1} + \alpha_{i,n-i})x_1 + \cdots +  (\alpha_{i,n-i-2} + \alpha_{i,n-i}) x_{n-i-2} + \alpha_{i,n-i-1} x_{n-i-1} \big)},   
 \end{align*}
 yielding that  
 \[
   \alpha_{i+1,j} =  \alpha_{i,j} +  \alpha_{i,n-i},\qquad j\in\{1,\ldots,n-i-2\},
    \qquad  \alpha_{i+1,n-i-1} =  \alpha_{i,n-i-1},\qquad 
    \beta_{i+1} = \beta_i\alpha_{i,n-i}.
 \]
In words, for each $n\geq 3$, each row of the matrix $A_n$ can be obtained from its previous row in the following way:
\begin{itemize}
  \item[(a)] the zero entries are always kept (if there are zero entries),
  \item[(b)] the last nonzero entry is replaced by 0, 
  \item[(c)] the last but one nonzero entry is unchanged,
  \item[(d)] the entries not mentioned in the previous steps are increased by the last nonzero value mentioned in step (b).
 \end{itemize}
The procedure above is repeated in a consecutive way.
For example, 
 \begin{align*}
  & A_3 = \begin{bmatrix}
          3 & 2 & 1 \\
          4 & 2 & 0 \\
        \end{bmatrix},
  \qquad 
  A_4=\begin{bmatrix}
        4 & 3 & 2 & 1 \\
        5 & 4 & 2 & 0 \\
        7 & 4 & 0 & 0 \\
      \end{bmatrix},\\
 & A_5= \begin{bmatrix}
          5 & 4 & 3 & 2 & 1 \\
          6 & 5 & 4 & 2 & 0 \\
          8 & 7 & 4 & 0 & 0 \\
          12 & 7 & 0 & 0 & 0 \\
        \end{bmatrix},
   \qquad 
   A_6 = \begin{bmatrix}
           6 & 5 & 4 & 3 & 2 & 1 \\
           7 & 6 & 5 & 4 & 2 & 0 \\
           9 & 8 & 7 & 4 & 0 & 0 \\
           13 & 12 & 7 & 0 & 0 & 0 \\
           20 & 12 & 0 & 0 & 0 & 0 \\
         \end{bmatrix}.
 \end{align*}

{\sl Step 9.}
By induction, we show that, for the matrix $A_n\in\RR_+^{(n-1)\times n}$ (defined in Step 8), it holds that 
 \begin{itemize}
  \item[(i)] its last ($(n-1)^{\mathrm{th}}$) row takes the form $(G_n, G_{n-1},0,\ldots,0)\in\RR^n$, 
  \item[(ii)] the last nonzero entry of its $i^{\mathrm{th}}$ row is $G_i$ for each $i\in\{1,\dots,n-1\}$,
 \end{itemize}
 where $G_n$, $n\in\NN$, are introduced and studied in Step 6. 

The statement holds in case of $n=3$, since, in Step 8, we have seen that the last row of $A_3$ is $(4,2,0)$ 
 and the last entry of its first row is $1$, and $G_1=1$, $G_2=2$ and $G_3=4$.
Now suppose that the statement holds for the matrices $A_3,\ldots,A_n$, and we show that it holds for the matrix $A_{n+1}$ as well.
Then, using that the first row $(n+1,n,\ldots,1)$ of $A_{n+1}$ coincides with 
 $(n+1,(A_n)_{1,1},\ldots,(A_n)_{1,n})$, by the procedure of obtaining the matrix $A_{n+1}$ from $A_n$ (described in Step 8),
 we have that the right-top $(n-1)\times n$ submatrix of $A_{n+1}$ is $A_n$, and hence 
 the $n\times(n+1)$ matrix $A_{n+1}$ takes the form
 \[
   A_{n+1}
    =\begin{bmatrix}
      n+1 & (A_n)_{1,1} & (A_n)_{1,2} & \ldots & (A_n)_{1,n-2} & (A_n)_{1,n-1} & (A_n)_{1,n} \\
      (A_{n+1})_{2,1} & (A_n)_{2,1} & (A_n)_{2,2} & \ldots & (A_n)_{2,n-2} & (A_n)_{2,n-1} & (A_n)_{2,n} \\
      \vdots & \vdots & \vdots & \ddots & \vdots & \vdots & \vdots \\
      (A_{n+1})_{n-1,1} & (A_n)_{n-1,1} & (A_n)_{n-1,2} & \ldots & (A_n)_{n-1,n-2} & (A_n)_{n-1,n-1} & (A_n)_{n-1,n} \\
      (A_{n+1})_{n,1} & (A_{n+1})_{n,2} & (A_{n+1})_{n,3} & \ldots & (A_{n+1})_{n,n-1}  & (A_{n+1})_{n,n}  & (A_{n+1})_{n,n+1} \\
     \end{bmatrix}.
  \]
By the induction hypothesis, we have  
 \[
  (A_n)_{n-1,1} = G_n, \qquad (A_n)_{n-1,2} = G_{n-1}, \qquad (A_n)_{n-1,j}= 0, \qquad j\in\{3,\ldots,n\}.
 \]
By the procedure of obtaining the matrix $A_{n+1}$ from $A_n$ (described in Step 8), we get that $(A_{n+1})_{n,2} = G_n$,
 and, using Step 6 and the induction hypothesis as well,  
 \begin{align*}
 (A_{n+1})_{n-1,1} 
  &= n+1 +G_1+G_2+\cdots + G_{n-2}
  = n+1 + (F_3-1) + (F_4-1)+ \cdots + (F_n-1)\\
  &= n+1 + F_3+\cdots + F_n - (n-2)
  = 3 +  F_3+\cdots + F_n \\
  &= 1 + F_1+ F_2+ F_3+\cdots + F_n
  = 1+ G_n 
  = 1 + F_{n+2}-1
  = F_{n+2}. 
 \end{align*} 
Therefore, using again the procedure of obtaining the matrix $A_{n+1}$ from $A_n$,
 we have that 
 \begin{align*}
   (A_{n+1})_{n,1}&= F_{n+2} + G_{n-1}
     =  F_{n+2} + (F_1+\cdots + F_{n-1})
    = (F_{n}+ F_{n+1}) + (F_1+\cdots + F_{n-1})\\
   & = F_1+\cdots + F_{n+1}
    = G_{n+1},
 \end{align*}
 as desired.
Since $(A_n)_{n-1,2}=G_{n-1}$ and $(A_n)_{n-1,j}=0$, $j\in\{3,\ldots,n\}$, by the procedure of obtaining the matrix $A_{n+1}$ from $A_n$,
 we have that $(A_{n+1})_{n,j}=0$, $j\in\{3,\ldots,n+1\}$.
This completes the proof. 

{\sl Step 10.}
We verify that, for each $n\geq 3$, $n\in\NN$, and $x_1,x_2\in\RR_+$, 
 \begin{align}\label{conjecture}
   I_{n-2}(x_1,x_2) = \frac{1}{G_1\cdots G_{n-2}}
             \ee^{- ( G_nx_1 + G_{n-1}x_2)},
 \end{align}
 where $G_n$, $n\in\NN$, are introduced and studied in Step 6. 

It follows by \eqref{help_exp_2} with $i=n-2$, Step 9 and the fact that $\beta_i=G_1G_2\cdots G_i$, $i\in\{1,\ldots,n-2\}$,
 which can be checked as follows.
Since $\beta_0=1$ and $\beta_{i+1} = \beta_i\alpha_{i,n-i}$, $i\in\{0,1,\ldots,n-3\}$ (obtained in Step 8), 
 we have  
 \begin{align*}
  \beta_1 = \beta_0\alpha_{0,n} = 1\cdot 1=G_1,\qquad \beta_2 = \beta_1\alpha_{1,n-1} = G_1\cdot G_2,
 \end{align*}
 since $\alpha_{1,n-1}$ is the last nonzero entry of the second row of $A_n$, and it equals to $G_2$ (due to Step 9).
Similarly, using that $\alpha_{i,n-i}$ is the last nonzero entry of the $(i+1)^{\mathrm{th}}$ row of $A_n$, and it equals to $G_{i+1}$,
 $i\in\{0,\ldots,n-3\}$ (due to Step 9), by induction, one can conclude the formula for $\beta_i$, $i\in\{0,\ldots,n-2\}$.

{\sl Step 11.}
By Step 7 and \eqref{conjecture}, we get 
 \begin{align*}
  P_n & = n! \int_{0}^{\infty} \int_{0}^{\infty}  I_{n-2}(x_1,x_2) \, \dd x_2\, \dd x_1
        = \frac{n!}{G_1\cdots G_{n-2}} \int_{0}^{\infty} \int_{0}^{\infty} \ee^{ - G_nx_1 - G_{n-1}x_2} \, \dd x_2\, \dd x_1\\
      & = \frac{n!}{G_1\cdots G_{n-1}} \int_{0}^{\infty}  \ee^{ - G_nx_1 }  \, \dd x_1
        = \frac{n!}{G_1\cdots G_{n}}.
 \end{align*}
By Step 6, this implies the statement. 
\proofend

\section{Pick-up sticks for Pareto  distribution}\label{Sec_Par}

In this section, we investigate a variant of the pick-up sticks problem considered in Section \ref{Sec_exp} 
 by replacing the exponential distribution with Pareto type I distribution.

Let $n\in\NN$, $\alpha>0$, and $X_1,\ldots,X_n$ be independent 
 and identically distributed random variables such that 
 the common distribution function takes the form
 \[
   \PP(X_1<x)
       = \begin{cases}
                   1 - x^{-\alpha} & \text{if $x>1$,}\\
                   0 & \text{if $x\leq 1$,}
               \end{cases}
 \]
 that is, $X_1$ has a Pareto type I distribution with scale parameter $1$ and shape parameter $\alpha$.

Recall that $X_{(1)},\ldots,X_{(n)}$ satisfying $X_{(1)}\leq \cdots \leq X_{(n)}$ denote the order statistics of  $X_1,\ldots,X_n$.
Let
 \begin{align}\label{R_def}
  R_k:=\frac{X_{(k)}}{X_{(k+1)}},\qquad k=1,\ldots,n-1.
 \end{align}

In the following Lemma \ref{Lem_Pareto}, we collect some known facts about Pareto type I distribution 
 that will be used in this section.
Part (ii) of Lemma \ref{Lem_Pareto} is about a representation of order statistics for Pareto type I distribution, and it is
  a special case of a general result that can be found, e.g., in Ahsanullah et al. \cite[Chapter 4]{AhsNevSha}. 
For completeness, we provide proofs as well.

\begin{Lem}\label{Lem_Pareto}
The following assertions hold:
 \begin{itemize}
  \item[(i)] $X_1\distre U^{-\frac{1}{\alpha}}$, where $U$ is a random variable with uniform distribution on the interval $(0,1)$.  
            The random variable $\ln(X_1)$ has an exponential distribution with parameter $\alpha$.
  \item[(ii)] 
     \begin{align}\label{Pareto_representation}
        (X_{(1)}, \ldots, X_{(n)}) 
       \distre \left( \exp\left(\frac{1}{\alpha} \sum_{j=1}^i \frac{Z_j}{n-j+1}\right)\right)_{i=1}^n,
     \end{align}
      where $Z_1,\ldots,Z_n$ are independent and exponentially distributed 
      random variables with parameter $1$. 
 \item[(iii)] The random variables $R_k$, $k=1,\ldots,n-1$, defined in \eqref{R_def} are independent, and $R_k$ has a Beta distribution with parameters
              $(n-k)\alpha$ and $1$, i.e., the density function of $R_k$ takes the form $(n-k)\alpha x^{(n-k)\alpha -1}\bbone_{[0,1]}(x)$, $x\in\RR$,
              equivalently,   
              \[
              \PP(R_k<x) = \begin{cases}
                              0 & \text{if $x\leq 0$,}\\
                              x^{(n-k)\alpha} & \text{if $x\in(0,1)$,}\\
                              1 & \text{if $x\geq 1$,}   
                            \end{cases}   
                            \qquad k=1,\ldots, n-1.
              \] 
              Further, we get that                 
              \begin{align}\label{help_pareto_4}
                 X_{(k)} = \left( \prod_{j=k}^{\ell-1}R_j \right) X_{(\ell)},\qquad k\leq \ell, \quad k,\ell\in\{ 1,\ldots,n-1\},
              \end{align}
              with the convention $\prod_{j=k}^{k-1}:=1$.
              In particular, for each $k\in\{ 1,\ldots,n-1\}$, $X_{(k)}$ can be obtained from the maximum $X_{(n)} = \max\{X_1,\ldots,X_n\}$ 
             by multiplying with some independent Beta distributed random variables.
 \end{itemize} 
\end{Lem}

\noindent{\bf Proof.}
(i). We have that $\PP(U^{-\frac{1}{\alpha}} <x) = 0$ for $x\leq 1$, and 
 \begin{align*}
    \PP(U^{-\frac{1}{\alpha}} <x)
       = \PP(U^{\frac{1}{\alpha}} > x^{-1})
       = \PP(U > x^{-\alpha})
       = 1- x^{-\alpha},\qquad x>1.
 \end{align*}
Further, we get
 \begin{align*}
 \PP(\ln(X_1) < y )
    = \PP(X_1 < \ee^{y} )
    = \begin{cases}
                   1 - \ee^{-\alpha y} & \text{if $y>0$,}\\
                   0 & \text{if $y\leq 0$,}
       \end{cases}
 \end{align*}   
  as desired. 
   
(ii).
Let $Y_i:=\ln(X_i)$, $i=1,\ldots,n$.
By part (i), the independence of $X_1,\ldots,X_n$ and the Borel measurability of the function $\ln$, we have that
 $Y_i$, $i=1,\ldots,n$, are independent and exponentially distributed random variables with parameter $\alpha$.
Since the function $\ln$ is strictly increasing, we get that $Y_{(j)} = \ln(X_{(j)})$, $j=1,\ldots,n$, i.e., $X_{(j)} = \ee^{Y_{(j)}}$, $j=1,\ldots,n$.
By R\'enyi's representation \eqref{Renyi_rep}, we get that
 \begin{align}\label{help_pareto_2}
   \big( Y_{(1)},\ldots, Y_{(n)}\big) 
     \distre \left( \frac{1}{\alpha} \sum_{j=1}^i \frac{Z_j}{n-j+1} \right)_{i=1}^n,
 \end{align}
 where $Z_1,\ldots,Z_n$ are independent, exponentially distributed random variables with parameter $1$.
Taking into account that $(X_{(1)},\ldots, X_{(n)}) = (\ee^{Y_{(1)}},\ldots, \ee^{Y_{(n)}})$
 and the Borel measurability of the mapping $\RR^n\ni (u_1,\ldots,u_n)\mapsto (\ee^{u_1},\ldots,\ee^{u_n})$, 
 \eqref{help_pareto_2} yields the statement of part (ii).

(iii). 
By part (ii) and the Borel measurability of the function $\ln$, we get
 \begin{align*}
   \big( \ln(X_{(1)}),\ldots, \ln(X_{(n)})\big) 
     \distre \left( \frac{1}{\alpha} \sum_{j=1}^i \frac{Z_j}{n-j+1} \right)_{i=1}^n,
 \end{align*}
 where $Z_1,\ldots,Z_n$ are independent, exponentially distributed random variables with parameter $1$.
Hence, by \eqref{R_def}, we have that
 \[
  \ln(R_k) = \ln(X_{(k)}) - \ln(X_{(k+1)}) \distre - \frac{Z_{k+1}}{(n-k)\alpha},\qquad k=1,\ldots,n-1,
 \] 
 which implies that 
 \[
  R_k \distre \ee^{-\frac{Z_{k+1}}{(n-k)\alpha}},\qquad k=1,\ldots,n-1.
 \]
Therefore, $R_1,\ldots,R_{n-1}$ are independent, and for each $k\in\{1,\ldots,n-1\}$, we get that
 $\PP(R_k<x)=0$, $x\leq 0$, $\PP(R_k<x)=1$, $x\geq 1$, and 
 \begin{align*}
  \PP(R_k<x)& = \PP\Big( \ee^{-\frac{Z_{k+1}}{(n-k)\alpha}} < x \Big)
             = \PP\left( -\frac{Z_{k+1}}{(n-k)\alpha} < \ln(x) \right)
             = \PP\left(Z_{k+1} > - (n-k)\alpha \ln(x)\right)\\
            & = \ee^{-(- (n-k)\alpha \ln(x))}
             = x^{(n-k)\alpha},
             \qquad x\in(0,1).
 \end{align*}
By \eqref{R_def}, we have that $X_{(k)} = R_k X_{(k+1)}$, $k=1,\ldots,n-1$, which, by induction, yields \eqref{help_pareto_4}.
\proofend

\begin{Thm}\label{thm4}
Let $n\in\NN$ be an integer greater than or equal to $3$.
If $n$ numbers are chosen independently according to a Pareto type I distribution with scale parameter $1$ and shape parameter $\alpha>0$
  from the interval $(1,\infty)$, then the probability that no three of the numbers can be the side lengths of a triangle is 
 \begin{align}\label{Pareto_formula_3}
 \begin{split}
 &(n-1)! \alpha^{n-1}\int_0^1\cdots\int_0^1
        \prod_{i=1}^{n-1} r_i^{(n-i)\alpha-1}
        \prod_{i=2}^{n-1} \bbone_{\big\{r_i\leq\frac1{1+r_{i-1}}\big\}}
         \,\dd r_{n-1}\,\cdots \dd r_1\\
 &\qquad =(n-1)! \alpha^{n-2} \int_0^1\cdots\int_0^1 \prod_{i=1}^{n-2} \left(s_i^{(n-i)\alpha-1} (D_i(s_1,\ldots,s_i))^{-\alpha}\right)  \,\dd s_1\,\cdots \dd s_{n-2},
 \end{split}
 \end{align}
 where $D_i(s_1,\ldots,s_i)$, $i=1,\ldots,n-2$, are given by the recursion $D_1(s_1):=1+s_1$, $D_2(s_1,s_2):=1+s_1+s_2$, and 
 \begin{align}\label{recursion_D}
 D_i(s_1,\ldots,s_i):=D_{i-1}(s_1,\ldots,s_{i-1}) + s_i D_{i-2}(s_1,\ldots,s_{i-2}),\qquad i=3,\ldots,n-2.
 \end{align} 
\end{Thm}

\noindent{\bf Proof of Theorem \ref{thm4}.}
Let $X_1,\ldots,X_n$ be independent Pareto type I distributed random variables with scale parameter $1$ and shape parameter $\alpha>0$.
Recall that $X_{(1)},\ldots,X_{(n)}$ satisfying $X_{(1)}\leq \cdots \leq X_{(n)}$ denote the order statistics of  $X_1,\ldots,X_n$.
By Steps 2 and 3 of the proof of Theorem \ref{thm2},
 the requirement that no three of $X_1, X_2, \ldots, X_n$  can be the side lengths of a triangle
 is equivalent to the condition that the inequality
 \begin{equation*}
    X_{(i)} + X_{(i+1)} \leq X_{(i+2)}
 \end{equation*}
 holds for all $i\in\{1,\ldots, n-2\}$.
Using \eqref{help_pareto_4}, these inequalities are equivalent to
 \begin{align*}
   \frac{X_{(i)}}{X_{(i+2)}} + \frac{X_{(i+1)}}{X_{(i+2)}} \leq 1
    \qquad \Longleftrightarrow \qquad 
      R_iR_{i+1} + R_{i+1} \leq 1
    \qquad \Longleftrightarrow \qquad   
    R_{i+1} \leq \frac{1}{1+R_i}
 \end{align*}
 for each $i\in\{1,\ldots, n-2\}$.
Consequently, the probability that no three of $X_1, X_2, \ldots, X_n$ can be the side lengths of a triangle is
 \begin{align*}
  \PP\left(\bigcap_{i=1}^{n-2} \left\{ R_{i+1} \leq \frac{1}{1+R_i} \right\} \right).
 \end{align*} 
Denoting the joint density function of $(R_1,\ldots,R_{n-1})$ by $f_{R_1,\ldots,R_{n-1}}$, and using its explicit form
 \begin{align*}
   f_{R_1,\ldots,R_{n-1}}(r_1,\ldots,r_{n-1}) 
        & = f_{R_1}(r_1)\cdots f_{R_{n-1}}(r_{n-1}) \\
        & = \alpha^{n-1}\prod_{i=1}^{n-1} (n-i) r_i^{(n-i)\alpha -1}\bbone_{[0,1]}(r_i),\qquad r_1,\ldots,r_{n-1}\in\RR,
 \end{align*}
 which follows from part (iii) of Lemma \ref{Lem_Pareto}, and using that 
 \[
   \bbone_{\left\{ \bigcap_{i=1}^{n-2} \left\{ R_{i+1} \leq \frac{1}{1+R_i} \right\} \right\} } 
      = \prod_{i=1}^{n-2} \bbone_{\left\{ R_{i+1} \leq \frac{1}{1+R_i} \right\}}
      =  \prod_{i=2}^{n-1} \bbone_{\big\{R_i\leq\frac1{1+R_{i-1}}\big\}},
 \]
 we have 
\begin{align}\label{help_pareto_6}
 \begin{split}
 &\PP\left(\bigcap_{i=1}^{n-2} \left\{ R_{i+1} \leq \frac{1}{1+R_i} \right\} \right)\\
 &\qquad  = \alpha^{n-1}\int_0^1\cdots\int_0^1
        \prod_{i=1}^{n-1} (n-i) r_i^{(n-i)\alpha-1}
        \prod_{i=2}^{n-1} \bbone_{\big\{r_i\leq\frac1{1+r_{i-1}}\big\}}
         \,\dd r_{n-1}\,\cdots \dd r_1\\
 &\qquad =(n-1)! \alpha^{n-1}\int_0^1\cdots\int_0^1
        \prod_{i=1}^{n-1} r_i^{(n-i)\alpha-1}
        \prod_{i=2}^{n-1} \bbone_{\big\{r_i\leq\frac1{1+r_{i-1}}\big\}}
         \,\dd r_{n-1}\,\cdots \dd r_1\\
 & \qquad =(n-1)! \alpha^{n-1}\lim_{\vare\downarrow 0}\int_\vare^1\cdots\int_\vare^1
        \prod_{i=1}^{n-1} r_i^{(n-i)\alpha-1}
        \prod_{i=2}^{n-1} \bbone_{\big\{r_i\leq\frac1{1+r_{i-1}}\big\}}
         \,\dd r_{n-1}\,\cdots \dd r_1. 
\end{split}         
 \end{align}
 
We will perform a change of variables in the multiple integral on the right-hand side of \eqref{help_pareto_6}.
For the sake of brevity, for all $(s_1,\ldots,s_{n-1})\in[0,1]^{n-1}$,  
let us introduce the notations $D_0:=1$, $D_1:=D_1(s_1)=1+s_1$, and $D_i:=D_i(s_1,\ldots,s_i)$, $i=2,\ldots,n-1$, 
 given by the recursion
\begin{align}\label{help_pareto_5}
 D_i:=D_{i-1} + s_i D_{i-2},\qquad i= 2,\ldots,n-1.
  \end{align} 
Note that $D_i\geq 1+s_i$, $i=1,\ldots,n-1$.
For completeness, we mention that in Theorem \ref{thm4}, we defined $D_i$ for $i=1,\ldots,n-2$, 
 and above we extended these definitions by introducing $D_0$ and $D_{n-1}$ as well.
 
Let $J_{n-1}:(-1,\infty)^{n-1}\to\RR^{n-1}$,
 \[
    J_{n-1}(u_1,\ldots,u_{n-1}):= (u_i(1+u_{i-1}))_{i=1}^{n-1},
        \qquad (u_1,\ldots,u_{n-1})\in(-1,\infty)^{n-1},
 \]  
 where $u_0:=0$.
Then $J_{n-1}$ is injective, 
 since if $J_{n-1}(u_1,\ldots,u_{n-1}) = J_{n-1}(\widetilde u_1,\ldots,\widetilde u_{n-1})$
 for some $(u_1,\ldots,u_{n-1}),(\widetilde u_1,\ldots,\widetilde u_{n-1})\in(-1,\infty)^{n-1}$, then 
 by solving recursively the system of equations
 \[
 u_i(1+u_{i-1}) = \widetilde u_i(1+\widetilde u_{i-1}), \qquad i=1,\ldots,n-1,
 \]
 we have that $u_i=\widetilde u_i$, $i=1,\ldots,n-1$.
For each $\vare\in[0,\frac{1}{2})$, let
 \[
   K^{(\vare)}_{n-1} :=\Big\{ (u_1,\ldots,u_{n-1}) \in[\vare,1]^{n-1} : u_i(1+u_{i-1})\leq 1, \;\; i=2,\ldots,n-1 \Big\}.
 \] 
We show that $J_{n-1}$ maps bijectively $K^{(\vare)}_{n-1}$ onto 
 \[
 J_{n-1}(K^{(\vare)}_{n-1}) = \left\{ (s_1,\ldots,s_{n-1})\in[0,1]^{n-1} : s_1\in[\vare,1], \; s_i\in\left[\vare\frac{D_{i-1}}{D_{i-2}},1\right],\;\; i=2,\ldots,n-1\right\}.
 \]
First, note that $\vare\frac{D_{i-1}}{D_{i-2}}\leq 1$ for each $i\in\{2,\ldots,n-1\}$ and $\vare\in[0,\frac{1}{2})$.
Indeed, if $i=2$, then $\frac{D_{i-1}}{D_{i-2}} = \frac{D_1}{D_0} = 1+ s_1 \leq 2$, and if $i\in\{3,\ldots,n-1\}$, then 
 \[
   \frac{D_{i-1}}{D_{i-2}} = \frac{D_{i-2} + s_{i-1} D_{i-3}}{D_{i-2}} = 1 + \frac{s_{i-1} D_{i-3}}{D_{i-2}}\leq 2,
 \]
 since $s_2D_0 = s_2\leq 1\leq 1+s_1 = D_1$ and 
 \[
 s_{i-1} D_{i-3} \leq D_{i-3} \leq D_{i-3} + s_{i-2}D_{i-4} = D_{i-2},\qquad i\in\{4,\ldots,n-1\}.
 \]
The injectivity of $J_{n-1}$ has already been checked.
Now we prove that  $J_{n-1}$ is surjective.
If $(s_1,\ldots,s_{n-1})\in[0,1]^{n-1}$ is such that
 $s_1\in[\vare,1]$ and $s_i\in[\vare\frac{D_{i-1}}{D_{i-2}},1]$, $i=2,\ldots,n-1$, then 
 the system of equations $s_i = u_i(1+u_{i-1})$, $i=1,\ldots,n-1$ (where $u_0:=0$) has exactly one solution 
 $(u_1,\ldots,u_{n-1})$ belonging to $K^{(\vare)}_{n-1}$, namely,
 \begin{align}\label{help_int_trans_1}
 u_1=s_1, \qquad  u_i  = \frac{s_i D_{i-2}}{D_{i-1}},
                     \qquad i=2,\ldots, n-1.
 \end{align}
Indeed, $u_1=u_1(1+0) = u_1(1+u_0)=s_1$, $u_2(1+u_1)=\frac{s_2D_0}{D_1}(1+s_1)=s_2$, and,
 by \eqref{help_pareto_5}, for $i\in\{3,\ldots,n-1\}$, we have that 
 \[
  u_i(1+u_{i-1}) = \frac{s_i D_{i-2}}{D_{i-1}}\left(1+ \frac{s_{i-1} D_{i-3}}{D_{i-2}} \right) 
      = \frac{s_i(D_{i-2} + s_{i-1} D_{i-3})}{D_{i-1}} = s_i,
 \]
 which shows that $(u_1,\ldots,u_{n-1})$ given by \eqref{help_int_trans_1} is indeed a solution of the system of equations
 $s_i=u_i(1+u_{i-1})$, $1,\cdots,n-1$ (with $u_0:=0$).
Since $J_{n-1}$ is injective, there does not exist another solution.
The solution in question belongs to  $K^{(\vare)}_{n-1}$, since
 $u_1=s_1\in[\vare,1]$ and 
 \[
  u_i = \frac{s_i D_{i-2}}{D_{i-1}}\geq \frac{\vare D_{i-1}}{D_{i-2}}\cdot \frac{D_{i-2}}{D_{i-1}} = \vare,\qquad i=2,\ldots,n-1.
 \]
By \eqref{help_int_trans_1}, we also have that
 \begin{align}\label{help_int_trans_6}
  J_{n-1}^{-1}(s_1,\ldots,s_{n-1}) = \left( s_1,\frac{s_2D_0}{D_1},\frac{s_3D_1}{D_2},\ldots, \frac{s_{n-1}D_{n-3}}{D_{n-2}} \right),
   \qquad (s_1,\ldots,s_{n-1})\in[0,1]^{n-1}.
 \end{align}
In particular, $J_{n-1}$ maps bijectively $K^{(0)}_{n-1}$ onto $[0,1]^{n-1}$, and $J_{n-1}(K^{(\vare)}_{n-1})$ is a 
 compact subset of $J_{n-1}((-1,\infty)^{n-1})$ for all $\vare\in[0,\frac{1}{2})$.
Hence the inverse $J_{n-1}^{-1}$ (which is defined on $J_{n-1}((-1,\infty)^{n-1})$) maps bijectively $J_{n-1}(K^{(\vare)}_{n-1})$ onto $K^{(\vare)}_{n-1}$
 for all $\vare\in[0,\frac{1}{2})$.
Further, $J_{n-1}$ is continuously differentiable, and, since the $i^{\mathrm{th}}$ coordinate of $J_{n-1}(u_1,\ldots,u_{n-1})$ is $u_i(1+u_{i-1})$,
 we have that the Jacobi matrix of $J_{n-1}$ at $(u_1,\ldots,u_{n-1})\in(-1,\infty)^{n-1}$ is lower triangular with diagonal 
 entries $1, 1+u_1,\ldots, 1+ u_{n-2}$.
Consequently, the determinant of the Jacobi matrix of $J_{n-1}$ at $(u_1,\ldots,u_{n-1})\in(-1,\infty)^{n-1}$ is 
 $\prod_{i=1}^{n-2} (1+u_i)\ne 0$ for all $(u_1,\ldots,u_{n-1})\in(-1,\infty)^{n-1}$.
By the inverse mapping theorem, $J_{n-1}((-1,\infty)^{n-1})$ is an open subset of $\RR^{n-1}$, 
 $J_{n-1}^{-1}$ is continuously differentiable and its Jacobi matrix has nonzero determinant 
 on $J_{n-1}((-1,\infty)^{n-1})$.
Further, by \eqref{help_int_trans_6}, we can also see that, for all $(s_1,\ldots,s_{n-1})\in[0,1]^{n-1}$,
 the $i^{\mathrm{th}}$ coordinate of $J_{n-1}^{-1}(s_1,\ldots,s_{n-1})$ 
 depends only on $s_1,\ldots,s_i$, $i=1,\ldots,n-1$,
 and hence the Jacobi matrix of $J_{n-1}^{-1}$ at  $(s_1,\ldots,s_{n-1})\in[0,1]^{n-1}$ 
 is lower triangular with diagonal entries $1$ and $\frac{D_{i-2}}{D_{i-1}}$, $i=2,\ldots,n-1$, yielding that  
 its determinant is equal to $\prod_{i=2}^{n-1} \frac{D_{i-2}}{D_{i-1}} = D_{n-2}^{-1}$.

Let us apply Jacobi's theorem (change of variables for multiple integrals) in Sagan \cite[Theorem 115.1, page 457]{Sag}
 for the integral
  \[
   (n-1)! \alpha^{n-1}\int_\vare^1\cdots\int_\vare^1
        \prod_{i=1}^{n-1} r_i^{(n-i)\alpha-1}
        \prod_{i=2}^{n-1} \bbone_{\big\{r_i\leq\frac1{1+r_{i-1}}\big\}}
         \,\dd r_{n-1}\,\cdots \dd r_1, \qquad \vare\in\left(0,\frac{1}{2}\right],
  \]
 with the choices $g:=J_{n-1}^{-1}$, $\Omega:=J_{n-1}((-1,\infty)^{n-1})$, $\cD^{(\vare)}:=J_{n-1}(K^{(\vare)}_{n-1})$ and
 $f:K^{(\vare)}_{n-1}\to \RR$, 
 \[
   f(r_1,\ldots,r_{n-1}):= (n-1)! \alpha^{n-1} \prod_{i=1}^{n-1} r_i^{(n-i)\alpha-1}
                           \prod_{i=2}^{n-1} \bbone_{\big\{r_i\leq\frac1{1+r_{i-1}}\big\}}
                          = (n-1)! \alpha^{n-1} \prod_{i=1}^{n-1} r_i^{(n-i)\alpha-1}
 \]
 for $(r_1,\ldots,r_{n-1})\in K^{(\vare)}_{n-1}$, where the second equality follows from the definition of $K^{(\vare)}_{n-1}$.
Then $f$ is bounded and continuous, and \eqref{help_pareto_6} and Jacobi's theorem yield that
 \begin{align}\label{help_int_trans_2}
 \begin{split}
 &\PP\left(\bigcap_{i=1}^{n-2} \left\{ R_{i+1} \leq \frac{1}{1+R_i} \right\} \right)\\
 &\qquad = \lim_{\vare\downarrow 0}\int_{K^{(\vare)}_{n-1}} f(r_1,\ldots,r_{n-1})\,\dd r_1\ldots\dd r_{n-1} \\
 &\qquad = \lim_{\vare\downarrow 0} \int_{g(\cD^{(\vare)})}  f(r_1,\ldots,r_{n-1})\,\dd r_1\ldots\dd r_{n-1} \\
 &\qquad = \lim_{\vare\downarrow 0} \int_{\cD^{(\vare)}} (f\circ g)(s_1,\ldots,s_{n-1}) \vert \det(\bA(s_1,\ldots,s_{n-1}))\vert \,\dd s_1\ldots\dd s_{n-1},
 \end{split}
 \end{align}
 where $\bA$ denotes the Jacobi matrix of $g=J_{n-1}^{-1}$.
As it was already explained, the determinant of the Jacobi matrix $\bA$ of $g=J_{n-1}^{-1}$ at a point $(s_1,\ldots,s_{n-1})\in \cD^{(\vare)}\subseteq [0,1]^{n-1}$ is
 \begin{align}\label{help_int_trans_3}
   \det(\bA(s_1,\ldots,s_{n-1})) = D_{n-2}^{-1}.
 \end{align}
Further, using \eqref{help_int_trans_6}, for all $(s_1,\ldots,s_{n-1})\in \cD^{(\vare)}$, we obtain 
 \begin{align*}
  (f\circ g)(s_1,\ldots,s_{n-1})
    & = f\left( s_1,\frac{s_2 D_0}{D_1}, \frac{s_3 D_1}{D_2},\ldots, \frac{s_{n-1} D_{n-3}}{D_{n-2}}  \right) \\
    & = (n-1)! \alpha^{n-1} s_1^{(n-1)\alpha -1} \prod_{i=2}^{n-1} \left( \frac{s_i D_{i-2}}{D_{i-1}}\right)^{(n-i)\alpha-1}\\
    & = (n-1)! \alpha^{n-1} \prod_{i=1}^{n-1} s_i^{(n-i)\alpha -1}
        \prod_{i=2}^{n-1} \left( \frac{D_{i-2}}{D_{i-1}}\right)^{(n-i)\alpha-1} \\
    & =  (n-1)! \alpha^{n-1} \prod_{i=1}^{n-1} s_i^{(n-i)\alpha -1} \left(\prod_{i=1}^{n-3} D_i^{-\alpha}\right)  D_{n-2}^{-(\alpha-1)}.
 \end{align*}
Hence, by \eqref{help_int_trans_2} and \eqref{help_int_trans_3}, we get 
 \begin{align}\label{help_int_trans_4}
  \begin{split}
  \PP\left(\bigcap_{i=1}^{n-2} \left\{ R_{i+1} \leq \frac{1}{1+R_i} \right\} \right)
   & = (n-1)! \alpha^{n-1} \lim_{\vare\downarrow 0}  \int_{\cD^{(\vare)}} \prod_{i=1}^{n-1} s_i^{(n-i)\alpha -1} \prod_{i=1}^{n-2} D_i^{-\alpha}
     \,\dd s_1\ldots\dd s_{n-1}\\
   & = (n-1)! \alpha^{n-1}\int_0^1\cdots\int_0^1
          \prod_{i=1}^{n-1} s_i^{(n-i)\alpha-1}
        \prod_{i=1}^{n-2} D_i^{-\alpha}\,\dd s_1 \,\cdots \dd s_{n-1},
  \end{split}        
 \end{align}
 where the second equality can be checked as follows.
The improper integral on the right hand side of \eqref{help_int_trans_4} is well-defined and finite, 
 since the inequalities $D_i\geq 1+s_i \geq 1$, $i\in\{1,\ldots,n-1\}$, yield that
 \begin{align}\label{help_int_trans_5}
  0 \leq \prod_{i=1}^{n-1} s_i^{(n-i)\alpha-1} \prod_{i=1}^{n-2} D_i^{-\alpha}
   \leq \prod_{i=1}^{n-1} s_i^{(n-i)\alpha-1}, \qquad (s_1,\ldots,s_{n-1})\in(0,1]^n,
 \end{align}
 and the improper integral $\int_0^1\cdots\int_0^1\prod_{i=1}^{n-1} s_i^{(n-i)\alpha-1} \,\dd s_1\,\cdots \dd s_{n-1}$
 is convergent (indeed, for each $i\in\{1,\ldots,n-1\}$, it holds that $\int_0^1 s_i^{(n-i)\alpha-1}\,\dd s_i = \frac{1}{(n-i)\alpha}$).
Using the non-negativity of the integrand (see the first inequality in \eqref{help_int_trans_5}), we have 
 that the value of the convergent improper integral on the right hand side of \eqref{help_int_trans_4} coincides with 
 \[
  (n-1)! \alpha^{n-1} \lim_{m\to\infty} \int_{L_m} \prod_{i=1}^{n-1} s_i^{(n-i)\alpha-1} 
                                         \prod_{i=1}^{n-2} D_i^{-\alpha}  \,\dd s_1\,\cdots \dd s_{n-1}
 \]
 for any sequence $L_m$, $m\in\NN$, of bounded simple sets of $(0,1]^{n-1}$ (i.e., 
 the boundary of $L_m$ has $(n-1)$-dimensional Lebesgue measure $0$) such that
 $L_i\subseteq L_j\subseteq (0,1]^{n-1}$ for $i\leq j$, $i,j\in\NN$, and 
 $\bigcup_{m=1}^\infty L_m = (0,1]^{n-1}$ (i.e., $L_m$, $m\in\NN$, is a so-called exhaustion of $(0,1]^{n-1}$),
 see, e.g., Proposition 5.2.5 in Georgiev and Zennir \cite{GeoZen}.  
Consequently, since $\cD^{(1/m)}$, $m\in\NN$, is an exhaustion of $(0,1]^{n-1}$, we have the second equality in 
 \eqref{help_int_trans_4}, as desired.
 
Finally, since $\prod_{i=1}^{n-2} D_i^{-\alpha}$ depends on $s_1,\ldots,s_{n-2}$, but not on $s_{n-1}$,
 by integrating with respect to $s_{n-1}$, we get 
 \begin{align*}
  \PP\left(\bigcap_{i=1}^{n-2} \left\{ R_{i+1} \leq \frac{1}{1+R_i} \right\} \right)
    = (n-1)! \alpha^{n-2}\int_0^1\cdots\int_0^1
           \prod_{i=1}^{n-2} (s_i^{(n-i)\alpha-1} D_i^{-\alpha}) \,\dd s_1 \,\cdots \dd s_{n-2},
 \end{align*}
 as desired.
\proofend

Next, we formulate a corollary of Theorem \ref{thm4} for $n=3$.

\begin{Cor}\label{cor}
If three numbers are chosen independently according to a Pareto type I distribution with scale parameter $1$ and shape parameter $\alpha>0$
 from the interval $(1,\infty)$,
 then the probability that they cannot be the side lengths of a triangle is
\begin{align}\label{Pareto_formula_1}
      2\alpha \int_0^1 \frac{x^{2\alpha-1}}{ (1+x)^\alpha} \,\dd x.
 \end{align}
In particular, if $\alpha$ is a positive integer, then the probability in question is equal to
 \begin{align}\label{Pareto_formula_2}
   2\alpha \left( \binom{2\alpha-1}{\alpha}(-1)^\alpha \ln(2) 
                   +\sum_{k\in\{0,\ldots,2\alpha-1\}\setminus\{\alpha\}}  \binom{2\alpha-1}{k} (-1)^k \frac{2^{\alpha-k}-1}{\alpha-k}   \right).
 \end{align}
 \end{Cor}

\noindent{\bf Proof.}
{\sl First proof of \eqref{Pareto_formula_1}.}
This proof is a direct application of Theorem \ref{thm4}, which is based on a representation of the ratios of successive order statistics for Pareto type I distribution
 in terms of some independent Beta distributions and the maximum of the observations (see \eqref{help_pareto_4}).
Using \eqref{Pareto_formula_3} with $n=3$, we get that the probability in question takes the form
  \begin{align*}
  2! \alpha \int_0^1 s_1^{2\alpha -1}  (D_1(s_1))^{-\alpha}\,\dd s_1
   =  2 \alpha \int_0^1 s_1^{2\alpha -1}  (1+s_1)^{-\alpha}\,\dd s_1.
 \end{align*}

\smallskip

{\sl Second proof of \eqref{Pareto_formula_1}.}
This proof is based on a representation of order statistics for Pareto distribution
 in terms of independent exponential distributions (see \eqref{Pareto_representation}).
Let $X_1$, $X_2$ and $X_3$ be independent Pareto type I distributed random variables with scale parameter $1$ and shape parameter $\alpha>0$.
Recall that their order statistics is denoted by $X_{(1)}, X_{(2)}, X_{(3)}$ satisfying $X_{(1)}\leq X_{(2)}\leq X_{(3)}$.
By Steps 2 and 3 in the proof of Theorem \ref{thm2},
 we need to calculate the probability
 \begin{align*}
 P_3:=\PP\big(X_{(1)} + X_{(2)} \leq X_{(3)} \big).           
 \end{align*}
Using part (ii) of Lemma \ref{Lem_Pareto}, we have that
 \begin{align}\label{help_pareto_1}
  (X_{(1)}, X_{(2)}, X_{(3)}) 
   \distre \left(  \exp\left(\frac{Z_1}{3\alpha}\right) , \exp\left(\frac{1}{\alpha}\left( \frac{Z_1}{3} + \frac{Z_2}{2}\right) \right), 
                   \exp\left(  \frac{1}{\alpha}\left( \frac{Z_1}{3} + \frac{Z_2}{2} + Z_3\right)    \right) \right),
 \end{align}
 where $Z_1$, $Z_2$, $Z_3$ are independent, exponentially distributed random variables with parameter $1$.

Let us introduce the random variables $\tZ_i:=\frac{Z_i}{4-i}$, $i=1,2,3$.
Then $\tZ_i$, $i=1,2,3$, are independent, and $\tZ_i$ has an exponential distribution with parameter $4-i$.
Using \eqref{help_pareto_1}, we have 
  \begin{align*}
  &\PP(X_{(1)} + X_{(2)}  \leq X_{(3)})\\
  &\qquad = \PP\left(\exp\left(\frac{1}{\alpha}\tZ_1 \right)
          + \exp\left(\frac{1}{\alpha}\left( \tZ_1 + \tZ_2 \right)  \right) 
        \leq \exp\left( \frac{1}{\alpha}\left( \tZ_1 + \tZ_2 + \tZ_3 \right)   \right) \right) \\
  &\qquad= \PP\left( 1 + \exp\left(\frac{1}{\alpha} \tZ_2 \right)
          \leq \exp\left(\frac{1}{\alpha} (\tZ_2 + \tZ_3)  \right) \right)
 \end{align*}         
 \begin{align*}         
  &\qquad=\PP\left( 1 + \exp\left( -\frac{1}{\alpha} \tZ_2 \right)
         \leq \exp\left(\frac{1}{\alpha} \tZ_3  \right) \right) \\
  &\qquad = \PP\left( \tZ_3 \geq \alpha \ln\left(1 + \ee^{-\frac{1}{\alpha}\tZ_2 }\right) \right).
\end{align*}
This formula also shows that for the calculation of the probability $\PP(X_{(1)} + X_{(2)}  \leq X_{(3)})$,
 interestingly, we only need two (and not three) independent exponentially distributed random variables.
Therefore, denoting the joint density function of $(\tZ_1,\tZ_2,\tZ_3)$ by $f_{\tZ_1,\tZ_2,\tZ_3}$ and using its explicit form
 \[
   f_{\tZ_1,\tZ_2,\tZ_3}(x_1,x_2,x_3) = f_{\tZ_1}(x_1)f_{\tZ_2}(x_2) f_{\tZ_3}(x_3)
                                      = 3!\ee^{-(3x_1 + 2x_2+ x_3)},\qquad x_1,x_2,x_3\in\RR_+,
 \] 
 we get 
 \begin{align*}
    P_3 &= 3!\int_{\big\{ (x_1,x_2,x_3)\in\RR_+^3 \; : \; x_3\geq \alpha \ln\left(1 + \ee^{-\frac{1}{\alpha}x_2 }\right) \big\} } 
               \ee^{-(3x_1 + 2x_2 + x_3)} \,\dd x_3\, \dd x_2\, \dd x_1\\
        &= 3! \int_{0}^{\infty} \int_{0}^{\infty} \int_{\alpha \ln\left(1 + \ee^{-\frac{1}{\alpha}x_2 }\right)}^{\infty} 
             \ee^{-(3x_1 + 2x_2 + x_3)}  \, \dd x_{3}\, \dd x_2\, \dd x_1\\
        &= 2 \int_{0}^{\infty} \int_{\alpha \ln\left(1 + \ee^{-\frac{1}{\alpha}x_2 }\right)}^{\infty} 
             \ee^{-( 2x_2 + x_3)}  \, \dd x_{3}\, \dd x_2 \\
      &= 2\int_0^\infty \ee^{-2x_2} \left(1 + \ee^{- \frac{1}{\alpha}x_2}\right)^{-\alpha}
        \,\dd x_2 
       = 2\int_0^\infty 
        \frac{1}{( \ee^{2x_2/\alpha} + \ee^{x_2/\alpha} )^\alpha}\,\dd x_2\\
      & = 2\int_1^\infty \frac{\alpha}{u(u^2+u)^\alpha}\,\dd u,          
 \end{align*}
 where, at the last equality, we used the substitution $\ee^{x_2/\alpha}=u$.
Hence, using the substitution $u=\frac{1}{x}$, we have   
 \begin{align*}
  P_3 = 2\int_1^\infty \frac{\alpha}{u (u^2(1 + u^{-1}))^\alpha}\,\dd u 
     = 2\alpha \int_1^\infty \frac{u^{-(2\alpha+1)}}{(1 + u^{-1})^\alpha}\,\dd u
     = 2\alpha \int_0^1 \frac{x^{2\alpha-1}}{(1+x)^\alpha} \,\dd x.
  \end{align*}   

Now we turn to the proof of \eqref{Pareto_formula_2}, i.e., the special case when $\alpha$ is a positive integer.
By the substitution $y = x+1$ and the binomial theorem, we have 
\begin{align*}
 P_3 & = 2\alpha \int_0^1 \frac{x^{2\alpha-1}}{(1+x)^\alpha}\,\dd x = 2\alpha \int_1^2 \frac{(y-1)^{2\alpha-1}}{y^\alpha}\,\dd y\\
     & = 2\alpha \int_1^2 \frac{1}{y^\alpha} \sum_{k=0}^{2\alpha-1} \binom{2\alpha-1}{k} y^{2\alpha-1-k} (-1)^k \,\dd y
        = 2\alpha \sum_{k=0}^{2\alpha-1} \binom{2\alpha-1}{k} (-1)^k \int_1^2 y^{\alpha-1-k}\,\dd y.
\end{align*} 
Hence, using that
 \[
   \int_1^2 y^{\alpha-k-1}\,\dd y 
     = \begin{cases}
         \frac{2^{\alpha-k}-1}{\alpha-k} & \text{if $k\ne \alpha$,}\\[1mm] 
         \ln(2) & \text{if $k = \alpha$,}
        \end{cases}  
 \]
 we obtain \eqref{Pareto_formula_2}.
\proofend

In the next remark, we point out that there are other approaches for determining the probability in Corollary \ref{cor}, 
 and we demonstrate it in case of $\alpha=1$ and $\alpha=\frac{1}{2}$.

\begin{Rem}\label{Rem1}
Let $X_1$, $X_2$ and $X_3$ be independent Pareto type I distributed random variables with scale parameter $1$ and shape parameter $\alpha>0$.
By part (i) of Lemma \ref{Lem_Pareto}, $X_i\distre U_i^{-\frac{1}{\alpha}}$, $i=1,2,3$, 
 where $U_i$, $i=1,2,3$, are independent random variable with uniform distribution on $(0,1)$.

(i). 
Suppose that $\alpha=1$.
Then $(X_1,X_2,X_3)\distre (U_1^{-1},U_2^{-1},U_3^{-1})$.
We need to calculate the probability
 \begin{align*}
 P_3&:=\PP\Big( \{ X_1 + X_2 \leq X_3\} \cup \{ X_1 + X_3 \leq X_2\} \cup \{ X_2 + X_3 \leq X_1\}  \Big)\\
    &= \PP\Big( \{ U_1^{-1} + U_2^{-1} \leq U_3^{-1}\} \cup \{ U_1^{-1} + U_3^{-1} \leq U_2^{-1}\} \cup \{ U_2^{-1} + U_3^{-1} \leq U_1^{-1}\}  \Big).
 \end{align*}
By the method of geometric probability, $P_3$ is equal to the $3$-dimensional Lebesgue measure of
 \[
   \Big\{ (x_1,x_2,x_3)\in(0,1)^3 : \{ x_1^{-1} + x_2^{-1} \leq x_3^{-1}\} \cup \{ x_1^{-1} + x_3^{-1} \leq x_2^{-1}\} 
                                      \cup \{ x_2^{-1} + x_3^{-1} \leq x_1^{-1}\}  \Big\}.
 \] 
By symmetry, we have that $P_3$ is equal to 3 times the $3$-dimensional Lebesgue measure of the set
 \[
   \Big\{ (x,y,z)\in(0,1)^3 :  x^{-1} + y^{-1} \leq z^{-1}  \Big\}
    = \Bigg\{ (x,y,z)\in(0,1)^3 : z\leq \frac{xy}{x+y} \Bigg\}.
 \]
Consequently, we have that 
 \begin{align*}
  P_3&=3\int_0^1\int_0^1 \frac{xy}{x+y} \,\dd y\,\dd x 
      = 3\int_0^1\int_0^1 \left( x - \frac{x^2}{x+y}\right)\,\dd y \,\dd x  
      = 3\int_0^1  \left[xy - x^2\ln(x+y) \right]_{y=0}^{y=1}\,\dd x  \\
     & = 3\int_0^1  \Big( x  - x^2\ln(x+1) + x^2\ln(x) \Big)\,\dd x  
       = 2(1-\ln(2)), 
  \end{align*}
 where the last equality follows by integration by parts.
Note that \eqref{Pareto_formula_2} for $\alpha=1$ also gives $2(1-\ln(2))$, as expected.

(ii). 
Suppose that $\alpha=\frac{1}{2}$.
Then $(X_1,X_2,X_3)\distre (U_1^{-2},U_2^{-2},U_3^{-2})$.
We need to calculate the probability
 \begin{align*}
 P_3&=\PP\Big( \{ X_1 + X_2 \leq X_3\} \cup \{ X_1 + X_3 \leq X_2\} \cup \{ X_2 + X_3 \leq X_1\}  \Big)\\
    &= \PP\Big( \{ U_1^{-2} + U_2^{-2} \leq U_3^{-2}\} \cup \{ U_1^{-2} + U_3^{-2} \leq U_2^{-2}\} \cup \{ U_2^{-2} + U_3^{-2} \leq U_1^{-2}\}  \Big).
 \end{align*}
By the method of geometric probability, $P_3$ is equal to the $3$-dimensional Lebesgue measure of
 \[
   \Big\{ (x_1,x_2,x_3)\in(0,1)^3 : \{ x_1^{-2} + x_2^{-2} \leq x_3^{-2}\} \cup \{ x_1^{-2} + x_3^{-2} \leq x_2^{-2}\} 
                                      \cup \{ x_2^{-2} + x_3^{-2} \leq x_1^{-2}\}  \Big\}.
 \] 
By symmetry, we have that $P_3$ is equal to 3 times
 the $3$-dimensional Lebesgue measure of the set
 \[
   \Big\{ (x,y,z)\in(0,1)^3 :  x^{-2} + y^{-2} \leq z^{-2}  \Big\}
    = \Bigg\{ (x,y,z)\in(0,1)^3 : z\leq \frac{xy}{\sqrt{x^2+y^2}} \Bigg\}.
 \]
Consequently, we have 
 \begin{align*}
  P_3&=3\int_0^1\int_0^1 \frac{xy}{\sqrt{x^2+y^2}} \,\dd y\,\dd x 
      = 3\int_0^1 \Big[ x\sqrt{x^2+y^2 }\Big]_{y=0}^{y=1}  \,\dd x 
      = 3\int_0^1 \left( x\sqrt{x^2+1} - x^2 \right)\,\dd x  \\
     &= 3\Big[ \frac{1}{3}(x^2+1)^{3/2} - \frac{x^3}{3} \Big]_0^1
     = 2(\sqrt{2}-1).
 \end{align*}
\proofend
\end{Rem}

In the next remark, we point out the fact that the probability \eqref{Pareto_formula_1} is strictly decreasing in $\alpha$.

\begin{Rem}
Let $X$ be a Pareto type I distributed random variable with scale parameter $1$ and shape parameter $\alpha>0$.
The larger $\alpha$, the greater the value of the distribution function of $X$ at a given point $x>1$. 
Roughly speaking, this means that the larger $\alpha$, the greater the probability that $X$ is close to 1.
Hence, given three independent random variables $X_1, X_2$ and $X_3$ having Pareto type I distribution 
with scale parameter $1$ and shape parameter $\alpha>0$, 
we expect that if we increase $\alpha$, then the probability that $X_1, X_2$ and $X_3$ are all close to 1 
is greater and greater, and thus the probability that they cannot be the side lengths of a triangle is less and less.
This heuristic argument and Corollary \ref{cor} suggest that the function 
 \[
    (0,\infty)\ni \alpha \mapsto  2\alpha \int_0^1 \frac{x^{2\alpha-1}}{(1+x)^\alpha}\,\dd x :=f(\alpha)
 \]
 should be strictly decreasing. 
Next, we check that it is indeed the case.
Using the substitution $x=\ee^{-u}$, we get that 
 \begin{align}\label{help_pareto_9}
 f(\alpha) =  2\alpha \int_0^\infty \frac{\ee^{-2\alpha u}}{(1+\ee^{-u})^\alpha}\,\dd u
           =  2\alpha \int_0^\infty\ee^{-\alpha g(u)}\,\dd u, \qquad \alpha>0,
 \end{align}
 where $g(u):=2u+\ln(1+\ee^{-u})$, $u\in\RR_+$.
We are going to prove that $f'(\alpha)<0$ for all $\alpha>0$. 
Using \eqref{help_pareto_9}, we have  
 \begin{align}\label{help_pareto_10}
 f'(\alpha) = 2\int_0^\infty \ee^{-\alpha g(u)}(1-\alpha g(u))\,\dd u,\qquad \alpha>0,
 \end{align}
 where the differentation with respect to $\alpha$ under the integral sign is permitted, since 
 \begin{align*}
  \int_0^\infty \vert \ee^{-\alpha g(u)} g(u) \vert \,\dd u
      \leq \int_0^\infty \ee^{-2\alpha u} (2u + \ln(1+\ee^{-u}))\,\dd u
      \leq \int_0^\infty (2u + \ln(2)) \ee^{-2\alpha u} \,\dd u
      <\infty, \quad \alpha>0.
 \end{align*}
Using that, for all $\alpha>0$, 
 \[
   \frac{\dd}{\dd u} (u\ee^{-\alpha g(u)}) = \ee^{-\alpha g(u)}(1-\alpha u g'(u)),\qquad u>0,
 \] 
 and $\lim_{u\to\infty}(u\ee^{-\alpha g(u)}) =0$, by the Newton--Leibniz theorem, we obtain 
 \[
   \int_0^\infty \ee^{-\alpha g(u)}(1-\alpha u g'(u))\, \dd u =0, \qquad \alpha>0.
 \]
Hence, using \eqref{help_pareto_10}, we get 
 \begin{align*}
   f'(\alpha) &= 2\int_0^\infty \ee^{-\alpha g(u)}(1-\alpha g(u))\,\dd u - 2\int_0^\infty \ee^{-\alpha g(u)}(1-\alpha u g'(u))\, \dd u\\
              &= 2\alpha \int_0^\infty \ee^{-\alpha g(u)}(ug'(u) - g(u))\,\dd u, \qquad \alpha>0.
 \end{align*}
Since 
 \[
  ug'(u) - g(u) = -\frac{u}{1+\ee^u} - \ln(1+\ee^{-u})< 0,\qquad u > 0,
 \]  
 we have  
 \[
 f'(\alpha) = -2\alpha \int_0^\infty \ee^{-\alpha g(u)}\left( \frac{u}{1+\ee^u}  + \ln(1+\ee^{-u}) \right)\,\dd u < 0, \qquad \alpha>0.
 \] 

Further, we have that
 \begin{align}\label{help_pareto_11}
  \lim_{\alpha\downarrow 0} f(\alpha) = 1 \qquad \text{and}\qquad 
  \lim_{\alpha\uparrow \infty} f(\alpha) = 0.
 \end{align}
Indeed, using \eqref{help_pareto_9} and the substitution $y=2\alpha u$, we get 
 \[
   f(\alpha) =  2\alpha \int_0^\infty (1+\ee^{-u})^{-\alpha} \ee^{-2\alpha u} \,\dd u
         = \int_0^\infty \ee^{-y} \left( 1 + \ee^{-\frac{y}{2\alpha}} \right)^{-\alpha} \,\dd y,
             \qquad \alpha>0,
 \]
 which implies \eqref{help_pareto_11} by the dominated convergence theorem. 
\proofend
\end{Rem}

\section*{Statements and Declarations}

The authors declare that they have no known competing financial interests
or personal relationships that could have appeared to influence the work presented in this paper. 
Data sharing is not applicable to this article as no datasets were generated or analyzed during the current study.


\begin{thebibliography}{99}

\bibitem{AhsNevSha}
\textsc{Ahsanullah, M., Nevzorov, V. B.} and \textsc{Shakil, M.} (2013).
\textit{An Introduction to Order Statistics}.
Atlantis Press, Paris, France.

\bibitem{BreCalLin}
\textsc{Brennan, M., Callow, N.} and \textsc{Lin, T. C.} (2026).
Fibonacci numbers and the probability of polygon formation using random length sticks.
\textit{arXiv 2604.27573}
 
\bibitem{GeoZen}
\textsc{Georgiev, S. G.} and \textsc{Zennir, K.} (2026).
\textit{Multiple Integrals in Calculus}.
De Gruyter, Berlin.

\bibitem{KonLkhTurUppHil}
\textsc{Kong, L., Lkhamsuren, L., Turner, A., Uppal, A.} and \textsc{Hildebrand, A. J.} (2013).
Random points, broken sticks, and triangles. Project report, \textit{Illinois Geometry Lab.}
URL: \url{https://studyres.com/doc/5452523/random-points--broken-sticks--and-triangles-project-report}
 
\bibitem{LiuWan}
\textsc{Liu, Y.} and \textsc{Wang, Z.} (2022).
Some generalizations of random broken and pick-up stick problems.
\textit{Journal of Physics: Conference Series}
\textbf{2287} Article 012003.
 
\bibitem{Muk}
\textsc{Mukerjee, R.} (2024).
A statistical approach to broken stick problems.
\textit{Communications in Statistics -- Theory and Methods}
\textbf{53(15)} 5629--5637.

\bibitem{Nev}
\textsc{Nevzorov, V. B.} (2001).
\textit{Records: Mathematical Theory}.
American Mathematical Society, Providence, Rhode Island.

\bibitem{PetTen}
\textsc{Petersen, T. K.} and \textsc{Tenner, B. E.} (2020).
Broken bricks and the pick-up sticks problem.
\textit{Mathematics Magazine}
\textbf{93(3)} 175--185.

\bibitem{Ren}
\textsc{R\'enyi, A.} (1953).
On the theory of order statistics.
\textit{Acta Mathematica Academiae Scientiarum Hungaricae}
\textbf{4} 191--231.

\bibitem{Sag}
\textsc{Sagan, H.} (1974).
\textit{Advanced Calculus of Real-Valued Functions of a Real Variable and Vector-Valued Functions of a Vector Variable}.
Houghton Mifflin Co.

\bibitem{ShoWel}
\textsc{Shorack, G. R.} and \textsc{Wellner, J. A.} (1986).
\textit{Empirical Processes with Applications to Statistics}.
John Wiley \& Sons, New York.

\bibitem{SudSunTreWan}
\textsc{Sudbury, A., Sun, A., Treeby, D.} and \textsc{Wang, E.} (2026+).
Pick-up sticks and the Fibonacci factorial.
To appear in \textit{The Fibonacci Quarterly}.
Available also at \textit{arXiv 2504.19911}

\bibitem{Ver1}
\textsc{Verreault, W.} (2022).
On the probability of forming polygons from a broken stick.
\textit{Statistics and Probability Letters}
\textbf{180} Article 109237.

\bibitem{Ver2}
\textsc{Verreault, W.} (2022).
MacMahon partition analysis: A discrete approach to broken stick problems.
\textit{Journal of Combinatorial Theory, Series A}
\textbf{187} Article 105571.
   
   
\end{thebibliography}
\end{document}